\documentclass[10pt]{article}
\usepackage{amsmath}
\usepackage{amsfonts}

\advance \textheight by \topskip
\advance \textheight by 1pt
\makeatother
\def\bea{\begin{eqnarray}}
\def\ena{\end{eqnarray}}

\def\lar{\longrightarrow}
\def\non{\nonumber}
\def\deg{\hbox{deg}}
\def\gr{\hbox{gr}}
\def\dim{\hbox{dim}}
\def\ch{\hbox{ch}}
\def\ker{\hbox{Ker}}
\def\bz{{\bar z}}
\def\bt{{\bar T}}
\def\sgn{\hbox{sgn}}

\newcommand{\bc}[2]{
\left(
\begin{array}{c}{#1}\\{#2}\end{array}
\right)}
\newcommand{\mapdown}[1]{\Big\downarrow
\rlap{$\vcenter{\hbox{$\scriptstyle#1\,$}}$}}
\newcommand{\qed}{\hbox{\rule[-2pt]{3pt}{6pt}}}
\newtheorem{prop}{Proposition}
\newtheorem{theorem}{Theorem}

\newtheorem{lemma}{Lemma}
\newtheorem{cor}{Corollary}

\title{
Differential Structure of Abelian Functions
}

\author{
Koji Cho\thanks{
e-mail: cho@math.kyushu-u.ac.jp} 
and Atsushi Nakayashiki\thanks{
e-mail: 6vertex@math.kyushu-u.ac.jp}\\
Department of Mathematics, Kyushu University\\
}

\date{}
\begin{document}
\maketitle
\begin{abstract}
The space of abelian functions of a principally polarized abelian
variety $(J,\Theta)$ is studied as a module over the ring ${\cal D}$
of global holomorphic differential operators on $J$. We construct
a ${\cal D}$ free resolution in case $\Theta$ is non-singular.
As an application, in the case of dimension 2 and 3,
we construct a new linear basis of the space of abelian
functions which are singular only on $\Theta$ in terms of logarithmic
derivatives of the higher dimensional $\sigma$-function.

\end{abstract}

\section{Introduction}
Let $(J,\Theta)$ be a $g$-dimensional principally polarized Abelian variety 
and $A$ the affine ring of $J-\Theta$. We express $J$ as the quotient
of the $g$-dimensional vector space by some lattice, $J={\mathbb C}^g/\Gamma$
and $\Theta$ as the zero locus of a theta function $\theta(z)$ with
$z=(z_1,...,z_g)$ being linear coordinates of ${\mathbb C}^g$.
Analytically $A$ is isomorphic to the ring of meromorphic functions on 
$J$ which have 
poles only on $\Theta$. Such functions can be considered as 
meromorphic and periodic functions on ${\mathbb C}^g$ which have
poles only on $(\theta(z)=0)$. Obviously if we differentiate such 
a function with respect to $z_i$ then we again get a function 
with the same property.
This means that $A$ becomes a module 
over the ring of differential operators 
${\cal D}={\mathbb C}[\partial_1,...,\partial_g]$, 
$\partial_i=\partial/\partial z_i$.
It is a very curious problem to determine generators and relations
of the ${\cal D}$-module $A$.
The aim of this paper is to study these problems for $(J,\Theta)$ with
$\Theta$ being non-singular.

The case of dimension one is known from 
the classical theory of elliptic functions.
In this case the structure of $A$ is very simple.
Let ${\wp}(z)$ be the Weierstrass elliptic function and $p$ the point
of the elliptic curve corresponding to $z=0$. 
In this case ${\cal D}={\mathbb C}[\partial]$, $\partial=\frac{d}{dz}$.
As a ${\cal D}$-module $A$ is generated by $1$ and ${\wp}(z)$.
More precisely $1$, ${\wp}(z)$, ${\wp}'(z)$, ${\wp}''(z)$, ...
give a ${\mathbb C}$-linear basis of $A$, 
where ${\wp}'(z)=\frac{d}{dz}{\wp}(z)$ etc.
This fact is incorporated in the beautiful addition formula
of Frobenius and Stickelberger:
\bea
&&
(-1)^{\frac{(n-1)(n-2)}{2}}\prod_{k=1}^{n-1}k!
\frac{\sigma(z_1+\cdots+z_n)\prod_{i<j}\sigma(z_i-z_j)}
{\prod_{j=1}^n\sigma(z_j)^n}
=
\left|
\begin{array}{ccc}
1&\cdots &1\\
{\wp}(z_1)&\cdots&{\wp}(z_n)\\
\vdots&\quad&\vdots\\
{\wp}^{(n-2)}(z_1)&\cdots&{\wp}^{(n-2)}(z_n)\\
\end{array}
\right|,
\label{FS}
\ena
where $\sigma(z)$ is the Weierstrass sigma function.
Consider both hand sides of (\ref{FS}) as a function of 
$z_1$. Then the left hand side is an element of
$A$ whose order of poles at $p$ is at most $n$.
The right hand side of (\ref{FS}) expresses the left hand
side as a linear combination of the basis $1$, ${\wp}'(z)$,..., 
${\wp}^{(n-2)}(z)$.

Up to now not many is known for the ${\cal D}$-module structure
of $A$ in the case of higher dimensions.
In \cite{NS1} the case of hyperelliptic Jacobians is studied and 
a conjecture on the ${\cal D}$-free resolution of $A$ is given.
Up to now it is still difficult to prove the conjecture in general.
Moreover few is known on the ${\cal D}$-module structure of $A$ for
non-hyperelliptic Jacobians.
In the present paper we construct a ${\cal D}$-free resolution of $A$
in the generic case which means that $\Theta$ is non-singular.
We remark that in this generic case similar problem for non-trivial
flat line bundles on $J$ is studied in \cite{N1}.

The content of the paper is as follows.
In section 2 the ${\cal D}$-module structure of the affine ring of
an abelian variety and its relation to the algebraic de Rham complex
are explained.  Large degree components of the highest cohomology 
group of the graded de Rham complex
are studied in section 3. As a consequence the affine ring $A$ 
is proved to be a finitely generated ${\cal D}$-module here.
In section 4 small degree components of the highest cohomology
group are studied.
The dimension of each homogeneous component of the highest cohomology
group are determined here.
In section 5 to 7 characters of cohomology groups, affine ring and
some related symplectic vector space are calculated.
A ${\cal D}$-free resolution is constructed in section 8.
In section 9 the results of \S 8 is interpreted into the term of 
theta functions. As examples a linear basis of $A$ is given in terms 
of logarithmic derivatives of a theta function in the case of genus 
two and three.
Three appendices provide proofs of Lemmas and assertions which are used
in the main body of the paper.

\section{Affine ring}
Let $(J,\Theta)$ be a principally polarized abelian variety 
of dimension $g$.
Throughout this paper we assume that $g\geq 2$ and $\Theta$ is non-singular.
This means, in particular, that Jacobian varieties of 
hyperelliptic curves of genus $g\geq 3$ and non-hyperelliptic curves
of genus $g\geq 4$ are excluded.

Let ${\cal O}$ be the sheaf of germs of holomorphic functions on $J$,
${\cal O}(n)$ $(n\geq 0)$ the sheaf of germs of meromorphic functions
on $J$ which have poles only on $\Theta$ of order at most $n$ and 
${\cal O}(\ast)$ the sheaf of germs of meromorphic functions
on $J$ which have poles only on $\Theta$. We set 
$A=H^0\left(J,{\cal O}(\ast)\right)$. It is isomorphic to the affine
coordinate ring of $J-\Theta$. The ring $A$ has an increasing filtration
determined by the order of poles on $\Theta$,
\bea
&&
A=\cup_{n=0}^\infty A_n,
\quad
A_n=H^0\left(J,{\cal O}(n)\right).
\non
\ena
We set $A_n=0$ for $n<0$ for convenience.

Analytically $A$ is decribed in the following manner.
Let $\tau$ be a $g$ by $g$ symmetric matrix whose imaginary part is 
positive definite, $J={\mathbb C}/{\mathbb Z}^g+\tau {\mathbb Z}^g$ and
$\theta(z)$ the Riemann's theta function,
\bea
&&
\theta(z)=
\sum_{n\in {\mathbb Z}^g}\exp\left(\pi i {}^tn\tau n+2\pi i{}^tn z\right),
\quad
z=(z_1,...,z_g).
\non
\ena
The theta divisor is defined by $\Theta=(\theta(z)=0)\subset J$.
Then
\bea
&&
A_n=\left\{\frac{f(z)}{\theta(z)^n}\right\},
\non
\ena
where $f(z)$ runs over all holomorphic functions on ${\mathbb C}^g$
with the property
\bea
&&
\frac{f(z+\tau p+q)}{\theta(z+\tau p+q)^n}=\frac{f(z)}{\theta(z)^n}
\label{periodic}
\ena
for any $p,q\in {\mathbb Z}^g$. Notice that the relation (\ref{periodic})
is preserved by the differentiation with respect to $z_i$. Thus $A$ becomes 
a module over the ring of differential operators 
${\cal D}={\mathbb C}[\partial_1,...,\partial_g]$ where 
$\partial_i=\partial/\partial z_i$. 

In order to study the ${\cal D}$-module structure of $A$ it is convenient
to consider the garded ring $\gr\,A$ associated with the filtration of $A$,
\bea
&&
\gr\,A=\oplus_{n=0}^\infty \gr_n\,A_n,
\quad
\gr_n\,A_n=A_n/A_{n-1}.
\non
\ena
Since the action of $\partial_i$ on $A$ satisfies the relation
\bea
&&
\partial_i A_n\subset A_{n+1},
\non
\ena
$\gr\,A$ also becomes a ${\cal D}$-module. 
Thus it is possible to define the deRham complex associated with $\gr\,A$
as follows. Let $T^\ast=\sum_{i=1}^g{\mathbb C}dz_i$ be the space
of translation invariant holomorphic one forms on $J$. Define the map
\bea
&&
d:\gr\,A\otimes \wedge^p T^\ast\lar \gr\,A\otimes \wedge^{p+1} T^\ast,
\non
\ena
by $d=\sum_{i=1}^g\partial_i\otimes dz_i$. It obviously defines a complex
$(\gr\,A\otimes \wedge^\bullet T^\ast,d)$. We denote by $H^k$ 
its $k$-th cohomology group $H^k\left(\gr\,A\otimes \wedge^\bullet T^\ast\right)$.
Notice that
\bea
&&
H^g\simeq \gr\,A/\sum_{i=1}^g\partial_i \gr\,A,
\non
\ena
as a vector space. A basis of this space gives a minimal set of generators
of $\gr\,A$ as a ${\cal D}$-module. Therefore we have to study
the cohomology groups of $(\gr\,A\otimes \wedge^\bullet T^\ast,d)$.
To this end we introduce a grading on $H^k$ as follows.
Let us assign degree to elements in $\wedge^p T^\ast$ by
\bea
&&
\deg\, dz_{i_1}\wedge\cdots \wedge dz_{i_g}=-p.
\non
\ena
The tensor product $\gr\,A\otimes \wedge^p T^\ast$ of two graded spaces is
naturally graded. Since the map $d$ preserves the grading, $H^k$ becomes graded. Let 
\bea
&&
H^k=\oplus_{n\geq -k}H^k_n
\non
\ena
be the decomposition into homogeneous comoponents.

\section{Cohomology of large degree}
We first study the case of $n$ large.
\begin{prop}\label{n-big1}
The following isomorphisms holds.
$$
H^g_{n}
\simeq
\left\{
\begin{array}{rl}
0,& \quad n\geq 2\\
H^{g}(J,{\cal O}),& \quad n=1.
\end{array}
\right.
$$
\end{prop}
\vskip2mm
\noindent

We set
\bea
&&
dz^g=dz_1\wedge\cdots\wedge dz_g.
\non
\ena
An element of $H^g_n$ can be written in the form $fdz^g$ for some 
$f\in \gr_{n+g}\, A$. We denote the vector space of such $f$'s by 
$H^g_n(dz^g)^{-1}$.

\begin{cor}\label{generated}
As a ${\cal D}$-module, $\gr\, A$ is generated by 
$\oplus_{n=-g}^{1}H^g_{n}(dz^g)^{-1}$.
\end{cor}

For the proof of Proposition \ref{n-big1} we need to introduce
some notations on sheaves.
Let $\Omega^p$ be the sheaf of germs of holomorphic $p$-forms
on $J$, $\Omega^p(n)$ $(n\geq 0)$ the sheaf of germs of meromorphic 
$p$-forms on $J$ which have poles only on $\Theta$ of order at most $n$,
$\Omega^p(-n)$ $(n\geq 0)$ the sheaf of germs of holomorphic 
$p$-forms on $J$ which have zeros on $\Theta$ of order at least $n$
and $\gr_n\,\Omega^p=\Omega^p(n)/\Omega^p(n-1)$.

Since $\Omega^p$ is a free ${\cal O}_J$-module, 
the following relations are valid
\bea
\gr_n\,\Omega^p&\simeq&\gr_n\,{\cal O}\otimes \Omega^p,
\label{free1}
\\
H^k\left(J,\gr_n\,\Omega^p\right)
&\simeq&
H^k\left(J,\gr_n\,{\cal O}\right)\otimes H^0\left(J,\Omega^p\right).
\label{free2}
\\
H^k\left(J,\Omega^p(n)\right)
&\simeq&
H^k\left(J,{\cal O}(n)\right)\otimes H^0\left(J,\Omega^p\right).
\label{free3}
\ena

To study cohomology groups of those sheaves the following
vanishing property of cohomologies due to Mumford \cite{M} is important.

\begin{lemma}\label{vanishing-lem}\cite{M}
We have
\bea
&&
H^k\left(J,{\cal O}(n)\right)=0,
\quad
k\geq 1, n\geq 1.
\label{vanishing}
\ena
\end{lemma}

\vskip5mm
The next lemma easily follows from this.
\vskip5mm

\begin{lemma}\label{fundamental-cohomology}
$$
H^i(J,\gr_{n}\,{\cal O})
\simeq
\left\{
\begin{array}{rl}
0,& \quad n<0,\quad i\leq g-2,
\\
H^i(J,{\cal O}),&
\quad n=0,\quad i\leq g-2,
\\
H^{i+1}(J,{\cal O}),&
\quad n=1,\quad i\geq 0,
\\
0,&
\quad n> 1,\quad i\geq 1
\end{array}
\right.
$$
\end{lemma}

Notice that $H^0\left(J,\Omega^p\right)\simeq 
\wedge^p H^0\left(J,\Omega^1\right)=\wedge^p T^\ast$.
Using (\ref{free2}), (\ref{free3}) and Lemma \ref{vanishing-lem}
we have
\bea
H^0\left(J,\Omega^p(n)\right)&\simeq& A_n\otimes \wedge^p T^\ast,
\label{further1}
\\
H^0\left(J,\gr_n\,\Omega^p\right)&\simeq& 
\gr_n\, A\otimes \wedge^p T^\ast, \quad n\geq 2.
\label{further2}
\ena

The exterior differentiation defines a
map $d:\Omega^p(n)\lar \Omega^{p+1}(n+1)$ which induces
a map $d:\gr_n\,\Omega^p\lar \gr_{n+1}\,\Omega^{p+1}$.
The induced map 
$d:H^0(J,\gr_n\,\Omega^p)\lar H^0(J,\gr_{n+1}\,\Omega^{p+1})$
on the cohomology groups
is the same as that of the complex $(\gr\,A\otimes \wedge^\bullet,d)$
due to the isomorphism (\ref{further2}).

Define the sheaf $\Phi^p_n$, $n\geq 1$ as the kernel of the map 
$d:\Omega^p(n)\lar \Omega^{p+1}(n+1)$.
By the definition the following sequence is exact
\bea
&&
0\lar \Phi^p_n \lar \gr_n\, \Omega^p \stackrel{d}{\lar}
 d\gr_n\,\Omega^p \lar 0.
\non
\ena

Notice that the map $d:\gr_n\,\Omega^p \lar \gr_{n+1}\,\Omega^{p+1}$ is
${\cal O}$ linear and the sheaf $\Phi^p_n$, $d\gr_n\,\Omega^p$ 
become a coherent ${\cal O}$-module.
Let us set $\Xi=d\log\,\theta$.
It defines a map
\bea
&&
\Xi\wedge: \gr_n\,\Omega^p \lar \gr_{n+1}\Omega^{p+1}.
\non
\ena

\begin{lemma}\label{sheaf-lemma}
(i) For $p\geq 1$
$$
\Phi^p_n\simeq 
\left\{
\begin{array}{rl}
d\gr_{n-1}\,\Omega^{p-1},
&
\quad 
n\geq 2,
\\
\Xi\wedge \gr_{0}\Omega^{p-1},
&
\quad
n=1
\end{array}
\right.
\non
$$
(ii) $\Xi\wedge \gr_{n}\,{\cal O}\simeq \gr_{n}\,{\cal O}$ for any integer $n$.
\vskip2mm

\noindent
(iii) $\gr_n{\cal O}\simeq d\gr_n{\cal O}$.
\vskip2mm

\noindent
(iv) 
$\ker(\Xi\wedge:\gr_{-n}\,\Omega^p\lar \gr_{-n+1}\,\Omega^{p+1})
\simeq \Xi\wedge \gr_{-n-1}\,\Omega^{p-1}$, $p\geq 1$, $n\geq 0$.
\end{lemma}
\noindent
{\it Proof.} Since every statement can be proved in a similar way, we
shall give a proof of (i).
It is obvious that the support of $\Phi^p_n$ is contained in $\Theta$.
Let $Q$ be a point of $\Theta$ and $(z_0,...,z_{g-1})$ be a local coordinate
system around $Q$ such that $z_0=0$ is a local defining equation of $\Theta$.
Write a local section of $\gr_n\,\Omega^p$ as
\bea
&&
\eta=\frac{1}{z_0^n}dz_0\wedge \eta_1+\frac{1}{z_0^n}\eta_2,
\non
\ena
where $\eta_1$, $\eta_2$ does not contain $dz_0$ and $z_0$.
Then, in $\gr_{n+1}\,\Omega^p$,
\bea
&&
d\eta=-\frac{n}{z_0^{n+1}}\wedge \eta_2.
\non
\ena
Thus $d\eta=0$ in $\gr_{n+1}\,\Omega^p$ is equivalent to
$\eta_2=0$. Therefore $\eta$ is a local section of $\Phi^p$
if and only if it is written as $\eta=\frac{1}{z_0^n}dz_0\wedge \eta_1$
with $\eta_1$ astisfying the condition above.
If $n\geq 2$, then 
$\eta=\frac{1}{1-n}d\left(\frac{\eta_1}{z_0^{n-1}}\right)$ in
$\gr_{n}\Omega^p$ and $\Phi^p_n=d\gr_{n-1}\,\Omega^{p-1}$.
If $n=1$, then $\eta=\frac{dz_0}{z_0}\wedge\eta_1$. 
Thus $\Phi^p_1=d\log \theta\wedge \gr_0\,\Omega^{p-1}$.
\hfill\qed

\vskip5mm

\noindent
{\it Proof of Proposition \ref{n-big1}.}
\par
\noindent
By the definition
\bea
&&
H^g_{n-g}=
\frac{\gr_n\, A\otimes \wedge^g T^\ast}
{d\left(\gr_{n-1}\,A \otimes \wedge^{g-1} T^\ast\right)}.
\label{def-hgn}
\ena
By (\ref{further2}) it is sufficient to prove
$$
\frac{H^0(J,\gr_n\,\Omega^g)}{dH^0(J,\gr_{n-1}\,\Omega^{g-1})}
\simeq
\left\{
\begin{array}{rl}
0,& \quad n\geq g+2\\
H^{g}(J,{\cal O}),& \quad n=g+1.
\end{array}
\right.
$$

Since $H^i(J,\gr_n\,\Omega^p)=0$ for $n\geq 2$, $i\geq 1$
by Lemma \ref{fundamental-cohomology} and (\ref{free2}),
the cohomology sequence of
\bea
&&
0 \lar d\gr_{n-2}\,\Omega^{g-2} \lar
\gr_{n-1}\,\Omega^{g-1} \stackrel{d}{\lar} \gr_n\,\Omega^g
\lar 0,
\quad
n\geq 3,
\label{seq1}
\ena
gives the isomorphism
\bea
&&
\frac{H^0(J,\gr_n\,\Omega^g)}{dH^0(J,\gr_{n-1}\,\Omega^{g-1})}
\simeq
H^1(J,d\gr_{n-2}\,\Omega^{g-2}),
\quad
n\geq 3.
\label{isom-1}
\ena
Using similar sheaf exact sequences we get, for $n\geq g+1$,
\bea
&&
H^1(J,d\gr_{n-2}\,\Omega^{g-2})
\simeq
H^2(J,d\gr_{n-3}\,\Omega^{g-3})
\simeq
\cdots
\simeq
H^{g-1}(J,d\gr_{n-g}\,{\cal O}).
\non
\ena
The last cohomology can be easily calculated as desired using 
Lemma \ref{sheaf-lemma} (iii) and Lemma \ref{fundamental-cohomology}.
\hfill \qed
\vskip5mm

\section{Cohomology of small degree}
Next we study $H^g_n$ for small $n$.
Let us set
\bea
&&
a^{(g)}_n=\dim\,H^g_{n-g}.
\non
\ena
Obviously $a^{(g)}_0=1$, $a^{(g)}_1=0$ and $a^{(g)}_2=2^g-1$. 
By Proposition \ref{n-big1} we have $a^{(g)}_{g+1}=1$ and 
$a^{(g)}_n=0$ for $n\geq g+2$. 
The remaining values of $a^{(g)}_n$ are given by 

\begin{prop}\label{agn}
For $2\leq n\leq g+1$ we have
\bea
&&
a^{(g)}_n=
(-1)^{n-1}\bc{g}{n-2}
+\sum_{i=2}^{n}(-1)^{n-i}\bc{g+1}{n-i} i^g
\non
\\
&&
+\sum_{i=0}^{n-3}(-1)^i
\Big(
\bc{g}{n-1}\bc{g}{n-2-i}-\bc{g}{n}\bc{g}{n-3-i}
\Big).
\label{agn-formula}
\ena
\end{prop}
\vskip2mm
\noindent
{\it Proof.}
The case of $n=2$ can be easily checked.
The proof for $n=g+1$ is given in Proposition \ref{agn-2} (i).
Therefore we assume $3\leq n\leq g$.

By (\ref{further2}) and (\ref{def-hgn})
\bea
&&
a^{(g)}_n=\dim\, H^0(J,\gr_{n}\,\Omega^{g})-
\dim\, dH^0(J,\gr_{n-1}\,\Omega^{g-1}).
\non
\ena
We have
\bea
&&
dH^0(J,\gr_{n-1}\,\Omega^{g-1})
\simeq
\frac{H^0(J,\gr_{n-1}\,\Omega^{g-1})}{H^0(J,d\gr_{n-2}\,\Omega^{g-2})},
\quad
n\geq 3,
\non
\ena
by the cohomology sequence of (\ref{seq1}).
The dimension of $A_n$ is known as ({\it cf.} \cite{M,M2})
\bea
&&
\dim\, A_n=n^g.
\quad
n\geq 1,
\label{dim-an}
\ena
Due to (\ref{further2}) we have
\bea
&&
\dim\, H^0(J,\gr_{n}\,\Omega^{g})=n^g-(n-1)^{g},
\quad
n\geq 2.
\label{dim-1}
\ena
Therefore
\bea
a^{(g)}_n&=&
\dim\, H^0\left(J,\gr_n\,\Omega^g\right)
-
\dim\, H^0\left(J,\gr_{n-1}\,\Omega^{g-1}\right)
+
\dim\, H^0\left(J,d\gr_{n-2}\,\Omega^{g-2}\right)
\non
\\
&=&
n^g-(g+1)(n-1)^g+g(n-2)^g
+
\dim\, H^0\left(J,d\gr_{n-2}\,\Omega^{g-2}\right),
\quad
n\geq 3.
\label{dim0}
\ena

Let us calculate the last term of (\ref{dim0}).

Consider the exact sequence
\bea
&&
0
\lar 
d\gr_{k-1}\,\Omega^{g-n+k-1}
\lar
\gr_{k}\,\Omega^{g-n+k}
\stackrel{d}{\lar}
d\gr_{k}\,\Omega^{g-n+k}
\lar
0,
\quad
k\geq 2.
\label{s-exact1}
\ena
The long cohomology sequence of (\ref{s-exact1}) gives the exact
sequence
\bea
&&
0
\lar
H^0\left(J,d\gr_{k-1}\,\Omega^{g-n+k-1}\right)
\lar
H^0\left(J,\gr_{k}\,\Omega^{g-n+k}\right)
\lar
H^0\left(J,d\gr_{k}\,\Omega^{g-n+k}\right)
\non
\\
&&
\lar
H^1\left(J,d\gr_{k-1}\,\Omega^{g-n+k-1}\right)
\lar
0,
\label{h-exact1}
\ena
and the isomorphisms
\bea
&&
H^i\left(J,d\gr_{k}\,\Omega^{g-n+k}\right)
\simeq
H^{i+1}\left(J,d\gr_{k-1}\,\Omega^{g-n+k-1}\right),
\quad
i\geq 1.
\label{isom1}
\ena
By (\ref{h-exact1}) we have
\bea
\dim\, H^0\left(J,d\gr_{k}\,\Omega^{g-n+k}\right)
&=&
-\dim\,H^0\left(J,d\gr_{k-1}\,\Omega^{g-n+k-1}\right)
+\dim\,H^0\left(J,\gr_{k}\,\Omega^{g-n+k}\right)
\non
\\
&&
+\dim\,H^1\left(J,d\gr_{k-1}\,\Omega^{g-n+k-1}\right).
\label{dim1}
\ena
Multiplying $(-1)^{n-k}$ to both hand sides of (\ref{dim1}) and taking
summation in $k$ from $2$ to $n-2$ we get
\bea
\dim\,H^0\left(J,d\gr_{n-2}\,\Omega^{g-2}\right)
&=&
(-1)^{n-1}\dim\,H^0\left(J,d\gr_{1}\,\Omega^{g-n+1}\right)
+
\sum_{k=2}^{n-2}(-1)^{n-k}\dim\, H^0\left(J,\gr_{k}\,\Omega^{g-n+k}\right)
\non
\\
&&
+
\sum_{k=2}^{n-2}(-1)^{n-k}\dim\,H^1\left(J,d\gr_{k-1}\,\Omega^{g-n+k-1}\right).
\label{dim2}
\ena
By (\ref{further2}) and (\ref{dim-1}) we know that
\bea
&&
\dim\, H^0\left(J,\gr_{k}\,\Omega^{g-n+k}\right)
=
\bc{g}{n-k}\left(k^g-(k-1)^g\right),
\quad
k\geq 2.
\label{dim3}
\ena
Let us determine the remaining part in (\ref{dim2}).

Using repeatedly the isomorphism (\ref{isom1}) we get
\bea
&&
H^1\left(J,d\gr_{k-1}\,\Omega^{g-n+k-1} \right)
\simeq
H^2\left(J,d\gr_{k-2}\,\Omega^{g-n+k-2} \right)
\simeq
\cdots
\simeq
H^{k-1}\left(J,d\gr_{1}\,\Omega^{g-n+1} \right),
\label{isom2}
\ena
for $k\geq 2$.
Therefore one has to calculate the dimension of 
$H^{i}\left(J,d\gr_{1}\,\Omega^{g-n+1}\right)$.

By Lemma \ref{sheaf-lemma} (i) the following sequence is exact
\bea
&&
0
\lar
\Xi \wedge \gr_0\, \Omega^{g-n}
\lar
\gr_1\, \Omega^{g-n+1}
\stackrel{d}{\lar}
d\gr_1\, \Omega^{g-n+1}
\lar
0.
\label{s-exact2}
\ena
The long cohomology exact sequence of (\ref{s-exact2}) is
\bea
&&
\cdots\lar
H^{i}\left(J, \Xi \wedge \gr_0\, \Omega^{g-n}\right)
\stackrel{\alpha}{\lar}
H^{i}\left(J,\gr_1\, \Omega^{g-n+1}\right)
\lar
H^{i}\left(J,d\gr_1\, \Omega^{g-n+1}\right)
\lar\cdots
\label{h-exact2}
\ena
Let us study the map $\alpha$.

\begin{lemma}\label{sublem1}
We have
\bea
&&
H^{i}\left(J, \Xi \wedge \gr_{-k}\, \Omega^{g-n-k}\right)=0,
\quad
0\leq i\leq n+k-2, \quad
1\leq k\leq g-n.
\non
\ena
\end{lemma}
\vskip2mm
\noindent
{\it Proof.}
We can assume $n<g$.
Let us prove the lemma by descending induction on $k$.
For $k=g-n$,
\bea
&&
H^{i}\left(J, \Xi \wedge \gr_{-(g-n)}\, {\cal O}\right)
\simeq
H^{i}\left(J,\gr_{-(g-n)}\, {\cal O}\right)
=0,
\quad
i\leq g-2,
\non
\ena
by Lemma \ref{sheaf-lemma} (i),(ii) and Lemma \ref{fundamental-cohomology}.
Suppose that the lemma holds from $k+1$ to $g-n$.

By Lemma \ref{sheaf-lemma} (iv) we have the exact sequence
\bea
&&
0
\lar
\Xi \wedge \gr_{-k-1}\, \Omega^{g-n-k-1}
\lar
\gr_{-k}\, \Omega^{g-n-k}
\stackrel{\Xi\wedge}{\lar}
\Xi \wedge \gr_{-k}\, \Omega^{g-n-k}
\lar
0,
\non
\label{s-exact3}
\ena
for $0\leq k\leq g-n-1$.
Then
\bea
&&
H^{i}\left(J, \Xi \wedge \gr_{-k}\, \Omega^{g-n-k}\right)=0,
\quad
0\leq i\leq n+k-2,
\non
\ena
by the cohomology sequenec of (\ref{s-exact3}), 
Lemma \ref{fundamental-cohomology} and the induction hypothesis.
Thus the lemma is proved.
\qed

\begin{lemma}\label{sublem2}
$H^{i}\left(J, \Xi \wedge \gr_{0}\, \Omega^{g-n}\right)
\simeq 
H^{i}\left(J,\gr_{0}\, \Omega^{g-n}\right)$, \quad
 $i\leq n-2$.
\end{lemma}
\vskip2mm
\noindent
{\it Proof.}
The lemma follows from the long cohomology sequence of (\ref{s-exact3}), $k=0$,
and Lemma \ref{sublem1}.
\qed

Let ${\bar T}^\ast=\oplus_{i=1}^g \mathbb{C}d\bz_i$, where $\bar{}$ denotes
the complex conjugation.
Then $H^i\left(J,{\cal O}\right)\simeq \wedge^i \bt^\ast$.
By Lemma \ref{fundamental-cohomology}, Lemma \ref{sublem2} and
(\ref{free3}) we have
\bea
H^{i}\left(J, \Xi \wedge \gr_{0}\, \Omega^{g-n}\right)
&\simeq&
\wedge^i \bt^\ast \wedge^{g-n}T^\ast,
\quad
i\leq n-2,
\label{hi-1}
\\
H^{i}\left(J,\gr_{1}\, \Omega^{g-n+1}\right)
&\simeq&
\wedge^{i+1} \bt^\ast \wedge^{g-n+1}T^\ast,
\quad
i\geq 0.
\label{hi-2}
\ena
Let 
\bea
&&
\hat{\omega}=
\pi \sum_{i,j=1}^g \left(Im(\tau)^{-1}\right)_{ij}d\bz_i\wedge dz_j
\in \bt^\ast \wedge T^\ast.
\non
\ena

\begin{lemma}\label{sublem3}
In the description of (\ref{hi-1}), (\ref{hi-2}) the map
$\alpha$ is given by the wedging $\hat{\omega}$.
\end{lemma}
\vskip2mm
\noindent
The proof of this lemma is given in Appendix A.
\vskip2mm

\noindent
{\bf Remark} From the cohomology sequence of
\bea
&&
0\lar \mathbb{C} \lar {\cal O}(1) \lar d{\cal O}(1) \lar 0,
\non
\ena
we have
\bea
&&
H^1(J,\mathbb{C})\simeq H^0\left(J,d{\cal O}(1)\right).
\label{h1c}
\ena
We set
\bea
&&
\zeta_i(z)=\partial_i\log\,\theta(z).
\label{zeta-i}
\ena
Then the one form $d\zeta_j$ is naturally an element of the right hand side
of (\ref{h1c}) and it can be considered as an element of 
$H^1(J,\mathbb{C})\simeq \bt^\ast\oplus T^\ast$.
Let 
\bea
&&
\omega=\sum_{j=1}^g d\zeta_j\wedge dz_j
\label{omega}
\ena
be the element
of $H^2(J,\mathbb{C})\simeq \wedge^2 H^1(J,\mathbb{C})$.
Then we have
\bea
&&
\omega=\hat{\omega}
\quad
\mbox{in $H^2(J,\mathbb{C})$}.
\label{omega=omega}
\ena
The proof of (\ref{omega=omega}) is given in Appendix B.
\vskip2mm

By the standard argument using the representation theory of
$sl_2$ we have that the map
\bea
&&
\hat{\omega}\wedge:
\oplus_{i+j=n}\wedge^i \bt^\ast \wedge^j T^\ast
\lar
\oplus_{i+j=n+2}\wedge^i \bt^\ast \wedge^j T^\ast
\non
\ena
is injective for $n\leq g-1$. Thus we have

\begin{lemma}\label{sublem4}
The map $\alpha$ is injective for $0\leq i\leq n-2$.
\end{lemma}
\vskip2mm

By Lemma \ref{sublem4}, (\ref{h-exact2}) splits into exact sequences
\bea\
&&
0
\lar
H^{i}\left(J, \Xi \wedge \gr_0\, \Omega^{g-n}\right)
\stackrel{\alpha}{\lar}
H^{i}\left(J,\gr_1\, \Omega^{g-n+1}\right)
\lar
H^{i}\left(J,d\gr_1\, \Omega^{g-n+1}\right)
\lar
0,
\non
\\
&&
\quad
0\leq i\leq n-3,
\label{h-exact3}
\\
&&
0
\lar
H^{n-2}\left(J, \Xi \wedge \gr_0\, \Omega^{g-n}\right)
\stackrel{\alpha}{\lar}
H^{n-2}\left(J,\gr_1\, \Omega^{g-n+1}\right)
\lar
H^{n-2}\left(J,d\gr_1\, \Omega^{g-n+1}\right)
\lar\cdots.
\label{h-exact4}
\ena
From (\ref{hi-1}), (\ref{hi-2}), (\ref{h-exact3}) we get
\bea
\dim\, H^i\left(J, d\gr_1\, \Omega^{g-n+1}\right)
&=&
\dim\, \wedge^{i+1} \bt^\ast \wedge^{g-n+1}T^\ast
-
\dim\, \wedge^i \bt^\ast \wedge^{g-n}T^\ast
\non
\\
&=&
\bc{g}{i+1}\bc{g}{n-1}-\bc{g}{i}\bc{g}{n},
\quad
0\leq i\leq n-3.
\label{dim4}
\ena
We substitute (\ref{dim3}), (\ref{dim4}), using (\ref{isom2}),
into (\ref{dim2}) we get
\bea
&&
\dim\, H^0\left(J,d\gr_{n-2}\,\Omega^{g-2}\right)
\non
\\
&&
=\sum_{k=1}^{n-2}(-1)^{n-k}\dim\,H^{k-1}\left(J, d\gr_1\,\Omega^{g-n+1}\right)
+
\sum_{k=2}^{n-2}(-1)^{n-k}\dim\,H^{0}\left(J, \gr_k\,\Omega^{g-n+k}\right)
\non
\\
&&
=\sum_{k=0}^{n-3}(-1)^{n-3-k}
\left(
\bc{g}{n-1}\bc{g}{k+1}-\bc{g}{n}\bc{g}{k}
\right)
+
\sum_{k=2}^{n-2}(-1)^{n-k}\bc{g}{n-k}\left(k^g-(k-1)^g\right)
\non
\\
&&
=\sum_{k=0}^{n-3}(-1)^{k}
\left(
\bc{g}{n-1}\bc{g}{n-2-k}-\bc{g}{n}\bc{g}{n-3-k}
\right)
+
\bc{g}{2}(n-2)^g
\non
\\
&&
+
(-1)^{n-1}\bc{g}{n-2}
+
\sum_{i=2}^{n-3}(-1)^{n-i}\bc{g+1}{n-i}i^g.
\label{dim5}
\ena
Then the proposition follows from (\ref{dim0}) and (\ref{dim5}).
\hfill\qed

\begin{prop}\label{agn-2}
(i) The right hand side of (\ref{agn-formula}) is equal to one.

\noindent
(ii) $$
\sum_{n=0}^{g+1}a^{(g)}_n
=
\bc{2g}{g}-\bc{2g}{g-2}+g!-
\frac{(2g)!}{g!(g+1)!}.
$$
\end{prop}
\vskip2mm

\noindent
{\it Proof}
(1) We have
\bea
a^{(g)}_{g+1}
&=&
(-1)^g\bc{g}{g-1}
+
\sum_{i=2}^{g+1}(-1)^{g+1-i}\bc{g+1}{g+1-i}i^g
+
\sum_{i=0}^{g-2}(-1)^i\bc{g}{g-1-i}
\non
\\
&=&
(-1)^{g+1}+\sum_{i=1}^{g+1}(-1)^{g+1-i}\bc{g+1}{i}i^g
+
\sum_{i=0}^{g-2}(-1)^i\bc{g}{i}
\non
\\
&=&1,
\non
\ena
where in the last equality we use
\bea
&&
\sum_{i=0}^g(-1)^i\bc{g}{i}=0,
\quad
0=\left(x\frac{d}{dx}\right)^g(1+x)^{g+1}\vert_{x=-1}=
\sum_{i=0}^{g+1}(-1)^i\bc{g+1}{i}i^g.
\non
\ena

\noindent
(2) 
Write
\bea
&&
\sum_{n=2}^{g+1}a^{(g)}_n=S_1+S_2+S_3,
\non
\\
&&
S_1=\sum_{n=2}^{g+1}(-1)^{n-1}\bc{g}{n-2},
\non
\\
&&
S_2=\sum_{n=2}^{g+1}\sum_{i=2}^{n}(-1)^{n-i}\bc{g+1}{n-i} i^g,
\non
\\
&&
S_3=\sum_{n=3}^{g+1}\sum_{i=0}^{n-3}(-1)^i
\Big(
\bc{g}{n-1}\bc{g}{n-2-i}-\bc{g}{n}\bc{g}{n-3-i}
\Big).
\non
\ena
The sum $S_1$ is easily evaluated as
\bea
&&
S_1=(-1)^g.
\non
\ena
To calculate $S_3$ we use
\bea
&&
\sum_{n-m=k}\bc{g}{n}\bc{g}{m}=\bc{2g}{g-k},
\non
\ena
which follows from
\bea
&&
(1+x)^g(1+x^{-1})^g=x^{-g}(1+x)^{2g}.
\non
\ena
We have
\bea
&&
S_3=-1+\bc{2g}{g-1}-\bc{2g}{g-2}.
\non
\ena
We use
\bea
&&
(1-x)^{g+1}\left(x\frac{d}{dx}\right)^g(1-x)^{-1}
=\sum_{i=1}^{g+1}\sum_{n=i}^{g+1}(-1)^{n-i}\bc{g+1}{n-i}i^gx^n.
\non
\ena
to deduce 
\bea
&&
S_2=g!+(-1)^{g+1}.
\non
\ena
The assertion follows from these results.
\hfill \qed
\vskip2mm

It is known that (\cite{NS1})
\bea
&&
\dim\, H^g\left(J-\Theta,\mathbb{C}\right)=
\bc{2g}{g}-\bc{2g}{g-2}+g!-
\frac{(2g)!}{g!(g+1)!}.
\non
\ena
By the algebraic de Rham theorem
\bea
&&
H^i(J-\Theta,\mathbb{C})\simeq H^i(A\otimes \wedge^{\cdot}T^\ast,d).
\non
\ena
Then the above proposition gives a relaion between
the highest cohomology group
of the algebraic de Rham complex $(A\otimes \wedge^{\cdot}T^\ast,d)$ and
that of the graded complex $(\gr\,A\otimes \wedge^{\cdot}T^\ast,d)$.

\begin{cor}
There is an isomorphism
\bea
&&
H^g
\simeq
H^g\left(A\otimes \wedge^{\cdot}T^\ast\right).
\non
\ena
\end{cor}
\vskip2mm
\noindent
{\it Proof.}
Take a set of representatives $\{\eta_i dz^g\}$, $\eta_i\in \gr_{n_i}\, A$, of
a basis of $H^g$.
Let $f_i\in A_{n_i}$ be a representative of $\eta_i$. Then
we have a surjective map
\bea
H^g
&\lar&
H^g\left(A\otimes \wedge^{\cdot}T^\bullet\right),
\non
\\
\eta_i dz^g &\mapsto& f_i dz^g.
\non
\ena
By Proposition $\ref{agn-2}$ this map is also injective.
\hfill\qed
\vskip2mm

We list some values of $a^{(g)}_n$ for small $g$:
\bea
&&
a^{(2)}_0=1,
\quad
a^{(2)}_1=0,
\quad
a^{(2)}_2=3,
\quad
a^{(3)}_3=1,
\non
\\
&&
\dim\, H^2(J-\Theta,\mathbf{C})=5=1+3+1,
\non
\\
&&
a^{(3)}_0=1,
\quad
a^{(3)}_1=0,
\quad
a^{(3)}_2=7,
\quad
a^{(3)}_3=6,
\quad
a^{(3)}_4=1,
\non
\\
&&
\dim\, H^3(J-\Theta,\mathbf{C})=15=1+7+6+1,
\non
\\
&&
a^{(4)}_0=1,
\quad
a^{(4)}_1=0,
\quad
a^{(4)}_2=15,
\quad
a^{(4)}_3=25,
\quad
a^{(4)}_4=10,
\quad
a^{(4)}_5=1,
\non
\\
&&
\dim\, H^4(J-\Theta,\mathbf{C})=52=1+15+25+10+1,
\non
\\
&&
a^{(5)}_0=1,
\quad
a^{(5)}_1=0,
\quad
a^{(5)}_2=31,
\quad
a^{(5)}_3=96,
\quad
a^{(5)}_4=66,
\quad
a^{(5)}_5=15,
\quad
a^{(5)}_6=1,
\non
\\
&&
\dim\, H^5(J-\Theta,\mathbf{C})=210=1+31+96+66+15+1.
\non
\ena

\section{The character of $H^g$}
In general the character of a graded vector space $H=\oplus_n H_n$ is defined
by
\bea
&&
\ch\,H=\sum_{n}t^n\dim\,H_n.
\non
\ena
By the definition the character of $H^g$ is 
\bea
&&
\ch\,H^g=\sum_{n=0}^{g+1} a^{(g)}_n t^{n-g}.
\non
\ena

\begin{prop}\label{ch-hg}
\bea
t^g\ch\,H^g&=&
1-(1-t)^g+(1-t)^{g+1}
\left(
1+(t\frac{d}{dt})^g(1-t)^{-1}
\right)
\non
\\
&&
-\sum_{I_1}(-1)^i\bc{g}{m+i}\bc{g}{m}t^{g+1-m}
+\sum_{I_2}(-1)^i\bc{g}{m+i+2}\bc{g}{m}t^{g-m},
\non
\ena
where
\bea
&&
I_1=\{(m,i)\vert \, m+i\leq g-1, 0\leq m, 1\leq i\},
\non
\\
&&
I_2=\{(m,i)\vert \, m+i\leq g-2, 0\leq m, 1\leq i\}.
\non
\ena
\end{prop}
\vskip2mm

\noindent
{\it Proof.}
By Proposition \ref{agn} we have
\bea
&&
t^g\ch\, H^g=1+S'_1+S'_2+S'_3,
\non
\\
&&
S'_1=\sum_{n=2}^{g+1}(-1)^{n-1}\bc{g}{n-2}t^n,
\non
\\
&&
S'_2=\sum_{n=2}^{g+1}\sum_{i=2}^{n}(-1)^{n-i}\bc{g+1}{n-i} i^gt^n,
\non
\\
&&
S'_3=\sum_{n=3}^{g+1}\sum_{i=0}^{n-3}(-1)^i
\Big(
\bc{g}{n-1}\bc{g}{n-2-i}-\bc{g}{n}\bc{g}{n-3-i}
\Big)t^n.
\non
\ena
It is easy to show that
\bea
&&
S'_1=-t^2(1-t)^g+(-1)^gt^{g+2}.
\non
\ena
By a similar calculation to that of $S_2$ in (ii) of 
Proposition \ref{agn-2} we get
\bea
&&
S'_2=(1-t)^{g+1}\left(t\frac{d}{dt}\right)^g(1-t)^{-1}
-t(1-t)^{g+1}+(-1)^{g+1}t^{g+2}.
\non
\ena
The sum $S'_3$ can be written as
\bea
&&
S'_3=-\sum_{I_1}(-1)^i\bc{g}{m+i}\bc{g}{m}t^{g+1-m}
+\sum_{I_2}(-1)^i\bc{g}{m+i+2}\bc{g}{m}t^{g-m}.
\non
\ena
The proposition follows from these formulae.
\hfill \qed

\section{The space $W^k$}
Let
\bea
&&
V=\oplus_{i=1}^g\mathbf{C}\alpha_i\oplus_{i=1}^g\mathbf{C}\beta_i
=V_{-}\oplus V_{+}
\non
\ena
be a vector space of dimension $2g$.
We make $\wedge^p V$, $p\geq1$ a graded vector space by specifying degree
as
\bea
&&
\deg(\beta_{i_1}\wedge\cdots\wedge \beta_{i_\ell}\wedge \alpha_{j_1}
\wedge\cdots\wedge \alpha_{j_k})=1-k,
\quad
\ell\geq 1,
\non
\\
&&
\deg(\alpha_{j_1}\wedge\cdots\wedge \alpha_{j_k})=-k.
\non
\ena
Let $\omega=\sum_{i=1}^g \beta_i\wedge \alpha_i$. 
Then $\deg\,\omega=0$.
Define
\bea
&&
W^k=\frac{\wedge^k V}{\omega \wedge^{k-2} V}.
\non
\ena
Since $\omega\wedge^{k-2} V$ is a graded subspace of $\wedge^k V$,
$W^k$ becomes a graded vector space. Its
character is given by

\begin{lemma}\label{ch-wk}
\bea
&&
\ch\, W^k=
\bc{g}{k}t^{-k}+\sum_{m=0}^{k-1}\bc{g}{k-m}\bc{g}{m}t^{1-m}
-\sum_{m=0}^{k-2}\bc{g}{k-2-m}\bc{g}{m}t^{-m}.
\non
\ena
\end{lemma}
\vskip2mm

\noindent
{\it Proof.}
Notice that
\bea
&&
\ch W^k=\ch\wedge^{k} V-\ch(\omega \wedge^{k-2} V).
\non
\ena
According as the decomposition
\bea
&&
\wedge^k V=
\wedge^k V_{-}\oplus
\oplus_{m=0}^{k-1}\wedge^{k-m}V_{+}\wedge^{m} V_{-}.
\non
\ena
we have
\bea
&&
\ch\wedge^{k} V=
\bc{g}{k}t^{-k}+\sum_{m=0}^{k-1}\bc{g}{k-m}\bc{g}{m}t^{1-m}.
\non
\ena
On the other hand we have at least one $d\zeta_i$ in
$\omega \wedge^{k-2} V$ and therefore
\bea
&&
\ch(\omega \wedge^{k-2} V)
=\sum_{m=0}^{k-2}\bc{g}{k-2-m}\bc{g}{m}t^{-m}.
\non
\ena
\hfill \qed

\section{Character of affine ring}
Let us calculate the character of $\gr\, A\otimes \wedge^g T^\ast$.
We have, by (\ref{dim-an}),
\bea
t^g\ch\left(\gr\,A\otimes \wedge^g T^\ast\right)
&=&\sum_{n=0}^\infty \dim\,( \gr_n\, A)t^n
\non
\\
&=&1+\sum_{n=2}^\infty\left(n^g-(n-1)^g\right)t^n
\non
\\
&=&(1-t)
\left(
1+(t\frac{d}{dt})^g(1-t)^{-1}
\right).
\non
\ena

We introduce a grading in ${\cal D}$ by assigning degree as
$\deg\,\partial_i=1$.
Then the character of ${\cal D}$ is given by
\bea
&&
\ch\,{\cal D}=(1-t)^{-g}.
\non
\ena
Set $H^i=W^i$ for $i<g$. 
At this point we encounter a remarkable identity.

\begin{cor}\label{chA=complex}
\bea
&&
\ch\left(\gr\,A\otimes \wedge^g T^\ast\right)
=\sum_{i=0}^g(-1)^i\ch\,( {\cal D}\otimes H^{g-i}).
\non
\ena
\end{cor}
\vskip2mm

\noindent
{\it Proof.}
By Lemma \ref{ch-wk} 
\bea
\sum_{i=1}^g(-1)^i\ch W^{g-i}
=&&t^{-g}\left((1-t)^g-1\right)
+\sum_{I_1}(-1)^i\bc{g}{m+i}\bc{g}{m}t^{1-m}
\non
\\
&&
-\sum_{I_2}(-1)^i\bc{g}{m+i+2}\bc{g}{m}t^{-m},
\non
\ena
where $I_1$,$I_2$ are the same as those in 
Proposition \ref{ch-hg}. Then the corollary
follows from Proposition \ref{ch-hg}.
\hfill \qed

\section{Free resolution of affine ring}
Define the map 
$d:{\cal D}\otimes\wedge^k  V \lar {\cal D}\otimes\wedge^{k+1} V$
by
\bea
&&
d(P\otimes \eta)=\sum_{i=1}^g \partial_iP\otimes \alpha_i\wedge \eta,
\non
\ena
and the map 
$\omega:{\cal D}\otimes\wedge^k  V \lar {\cal D}\otimes\wedge^{k+2} V$
by
\bea
&&
\omega(P\otimes \eta)=P\otimes \omega\wedge \eta.
\non
\ena
Obviously these two maps are ${\cal D}$-linear and commute. 
Thus $d$ induces a ${\cal D}$-linear map
$d:{\cal D}\otimes W^k\lar {\cal D}\otimes W^{k+1}$
and defines a complex $({\cal D}\otimes W^{\bullet},d)$:
\bea
&&
0\lar {\cal D} \lar {\cal D}\otimes W^1 
\lar \cdots \lar
{\cal D}\otimes W^{g}
\lar
0. 
\non
\ena

The following proposition has been proved in \cite{NS1}.

\begin{prop}\cite{NS1}\label{NS-W}
The complex $({\cal D}\otimes W^{\bullet},d)$ is exact at 
${\cal D}\otimes W^k$, $k\neq g$.
\end{prop}

\noindent
{\bf Remark} It can be easily checked that the map $d$ preserves degree. 
Notice that the degree in this paper and that in \cite{NS1} are different. 
\vskip2mm

Let $W^k=\oplus_n W^k_n$ be the homogeneous
decomposition.
Notice that 
\bea
&&
W^g_{n-g}=\frac{\wedge^{n-1}V_{+}\wedge^{g-n+1}V_{-}}
{\omega \wedge^{n-2}V_{+}\wedge^{g-n}V_{-}}\quad n>1,
\quad
W^g_{1-g}=0,
\quad
W^g_{-g}=\wedge^g V_{-}.
\label{wh-1}
\non
\ena
For $J=(j_1,\ldots,j_n)$ we set $|J|=n$ and $
dz_J=dz_{j_1}\wedge\cdots\wedge dz_{j_n}$ ect.

\begin{lemma}\label{zeta-pole}
The meromorphic $g$-form
$d\zeta_I\wedge dz_J$, $|I|=n-1\geq 1$, $|J|=g-n+1$ has a pole on $\Theta$ of order at most $n$ and thus can
be considered as an element of $\gr_n\, A\otimes \wedge^g T^\ast$.
\end{lemma}
\vskip2mm
\noindent
{\it Proof.} Let $I=(i_1,...,i_{n-1})$ and 
$J^c=\{1,2,...,g\}\backslash J=\{k_1,...,k_{n-1}\}$, $k_1<\cdots<k_{n-1}$.
Then
\bea
&&
d\zeta_I\wedge dz_J=c\det(\zeta_{i_pk_q})_{1\leq p,q\leq n-1}dz^g,
\label{det-expr}
\ena
where $\zeta_{ij}=\partial_i\zeta_j=\partial_i\partial_j\log\,\theta(z)$ and
$c=\pm1$.
We set $\theta_{i_1...i_l}(z)=\partial_{i_1}\cdots\partial_{i_l}\theta(z)$.
Substitutes the expression
\bea
&&
\zeta_{ij}=-\frac{\theta_i(z)\theta_j(z)}{\theta(z)^2}+
\frac{\theta_{ij}(z)}{\theta(z)},
\non
\ena
into the determinant in (\ref{det-expr}). 
Then one easily see that $d\zeta_I\wedge dz_J$ has poles on $\Theta$
of order at most $n$.\rightline{\qed}
\vskip3mm

By Lemma \ref{zeta-pole} we have a natural map
\bea
\wedge^{n-1} V_{+}\wedge^{g-n+1} V_{-} &\lar& H^g_{n-g},
\label{wh-2}
\\
\beta_I\wedge \alpha_J &\mapsto& d\zeta_I\wedge dz_J.
\non
\ena
Since $\omega=\sum_{i=1}^g \beta_i\wedge \alpha_i\in V_{+}\wedge V_{-}$
is mapped to
\bea
&&
\sum_{i,j=1}^g\partial_j\partial_i\log\,\theta\,dz_j\wedge dz_i=0,
\non
\ena
the map (\ref{wh-2}) induces a map
\bea
&&
\varphi_n:W^g_{n-g}\lar H^g_{n-g}.
\label{wh-3}
\ena

\begin{lemma}\label{injective}
The map $\varphi_n$ is injective for $n\geq 0$.
\end{lemma}
\vskip2mm

The proof of this lemma is given in the appendix.
\vskip2mm

In the following we sometimes identify $W^g$ 
as a graded subspace of $H^g$ by identifying
$\beta_I\wedge \alpha_J$ with $d\zeta_I\wedge dz_J$.
Let us fix a decomposition
\bea
&&
H^g=W^g\oplus U^g,
\non
\ena
where $U^g=\oplus_nU^g_n$ is a graded subspace of $H^g$. 
We fix a homogeneous basis $\{u_\beta\vert \beta\in B\}$ of $U^g$.
With the help of this basis we define the map
\bea
&&
ev:{\cal D}\otimes H^g \lar \gr\,A\otimes \wedge^g T^\ast,
\label{map-ev}
\ena
in the following way. 

Let ${\cal D}=\oplus_n{\cal D}_n$ be the homogeneous decomposition.
For $\alpha_I\wedge \beta_J\in \wedge^g V$, it is possible to
write $dz_I\wedge d\zeta_J=fdz^g$ with $f\in A_m$ where
$m=|J|+1$ if $J\neq\phi$, $m=0$ otherwise by Lemma \ref{zeta-pole}.
Define first a map $ev':\wedge^g V\lar \gr\,A\otimes \wedge^g T^\ast$ by
\bea
&&
ev'\left(P\otimes (\alpha_I\wedge \beta_J)\right)=P(f)\,dz^g\in 
\gr_{m+n}\,A\, dz^g,\quad P\in {\cal D}_n.
\non
\ena
Since $ev'(\omega\wedge^{g-2} V)=0$, it induces a well-defined 
map $ev:\,{\cal D}\otimes W^g
\lar \gr\, A \,dz^g$.

Next take $P\otimes u_\beta\in {\cal D}_n\otimes U^g_{m-g}$.
Write $u_\beta=f_\beta dz^g$ with $f_\beta\in gr_m\,A$.
Then we set
\bea
&&
ev(P\otimes u_\beta)=P(f_\beta)dz^g\in \gr_{m+n}\, A\, dz^g.
\non
\ena

Defining a map $d:{\cal D}\otimes W^{g-1}\lar {\cal D}\otimes H^g$ by
\bea
&&
dw=(dw,0)\in ({\cal D}\otimes W^g)\oplus ({\cal D}\otimes H^g),
\non
\ena
where $dw$ in the right hand side is the one defined previously,
we get the complex
\bea
&&
0 \lar {\cal D} \stackrel{d}{\lar} {\cal D}\otimes W^1 
\stackrel{d}{\lar}\cdots \lar {\cal D}\otimes W^{g-1} 
\stackrel{d}{\lar} {\cal D}\otimes H^g
\stackrel{ev}{\lar} \gr\, A \otimes \wedge^g T^\ast \lar 0.
\label{resolution}
\ena

\begin{theorem}\label{resol-grA}
The complex (\ref{resolution}) is exact and gives a ${\cal D}$-free
resolution of $\gr\, A \otimes \wedge^g T^\ast$.
\end{theorem}
\vskip2mm
\noindent
{\it Proof.}
By Proposition \ref{NS-W} the complex $({\cal D}\otimes W^{\bullet},d)$ 
is exact at ${\cal D}\otimes W^{k}$, $0\leq k\leq g-1$. 
Corollary \ref{generated} shows that
$ev$ is surjective.  By Corollary \ref{chA=complex}
\bea
&&
\ch\, (\gr\, A dz^g)=
\ch\,\left(Coker({\cal D}\otimes W^{g-1}\lar {\cal D}\otimes H^g)\right).
\non
\ena
It shows that the complex is exact at ${\cal D}\otimes H^g$. 
This completes the proof.
\hfill \qed
\vskip2mm

Simply omitting $dz^g$ from $\gr\, A \otimes \wedge^g T^\ast$
we get a ${\cal D}$-free resolution of $\gr\, A$ from (\ref{resolution}).
It is possible to give
a free resolution of $A$ itself using the above result.

Take a representatives of $\{f_\beta dz^g\vert \beta\in B_{n-g}\}$
from $A_ndz^g$ and denote them by the same letters.
The element $d\zeta_I\wedge dz_J$ is considered as an
element from $Adz^g$.
We define the map
\bea
&&
\tilde{ev}:{\cal D}\otimes H^g \lar A,
\non
\ena
simply by sending $P\otimes f\,dz^g$ to $P(f)\in A$.
Again we have the complex
\bea
&&
0 \lar {\cal D} \stackrel{d}{\lar} {\cal D}\otimes W^1 
\stackrel{d}{\lar}\cdots \lar {\cal D}\otimes W^{g-1} 
\stackrel{d}{\lar} {\cal D}\otimes H^g
\stackrel{\tilde{ev}}{\lar} A \lar 0,
\label{resolution2}
\ena
where the maps other than $\tilde{ev}$ are the same as before.

\begin{cor}\label{resol-A}
The complex (\ref{resolution2}) is exact and it gives a 
${\cal D}$-free resolution of $A$.
\end{cor}
\vskip2mm

\noindent
{\it Proof.}
The surjectivity of the map $\tilde{ev}$ follows from that
of $ev$. We have to prove $\ker\,\tilde{ev}=d({\cal D}\otimes W^{g-1})$.
Suppose that $v\in \ker\,\tilde{ev}$. Decompose $v$ into
homogeneous components as $v=v_n+v_{n-1}+\cdots$ with $\deg\, v_i=i$.
Since $\tilde{ev}(v)=0$, we have $\tilde{ev}(v_n)=-\tilde{ev}(v_{n-1})-\cdots
\in A_{n-1+g}$. Thus $ev(v_n)=0$. Since (\ref{resolution}) is exact, 
$v_n=dw_n$ in $\gr\,Adz^g$ for some $w_n\in ({\cal D}\otimes W^{g-1})_n$.
Then $\tilde{ev}(v-dw_n)=\tilde{ev}(v)=0$. Now $v-dw_n=v_{n-1}+v_{n-2}+\cdots$.
Thus repeating the same argument we finally have $v=dw$ for some
$w\in {\cal D}\otimes W^{g-1}$.
\hfill \qed

\section{Theta functions}
In this section we interpret the results of the previous section
in terms of theta and abelian functions.
For $I=(i_1,\ldots,i_r)$, $J=(j_1,\ldots,j_r)$ we define
the element $(I;J)$ of $A$ by
\bea
&&
(I;J)=\det(\zeta_{i_kj_l})_{1\leq k,l\leq r}.
\non
\ena
In particular $(i;j)=\zeta_{ij}$.
By the definition $(I;J)$ is anti-symmetric with respect to
the permutations of the indices of $I$ and $J$ respectively and
satisfy $(I;J)=(J;I)$.
These functions arise from the relation
\bea
&&
d\zeta_I\wedge dz_{J^c}=\sgn(J,J^c)(I;J)dz^g,
\non
\ena
where $J^c=(j_{r+1},...,j_g)$, $\{j_{r+1},...,j_g\}=\{1,2...,g\}\backslash J$
and $\sgn(J,J^c)$ is the sign of the permutation 
$(J,J^c): (1,2,...,g)\mapsto (J,J^c)$.
Let $\tilde{ev}'$ be the composition 
$\wedge^g V\lar W^g \stackrel{\tilde{ev}}{\lar} A$, which maps 
$\beta_I\wedge \alpha_{J^c}$ to $\sgn(J,J^c)(I;J)$.
Then ${\tilde ev}'\left(\omega\wedge^{g-2} V\right)=0$.
It is equivalent to the relations
\bea
&&
\sum_{k=1}^{r+2}(-1)^{k-1}(i_1,\ldots,i_r,j_k;j_1,\ldots,\hat{j_k},\ldots,j_{r+2})
=0.
\label{lin-rel}
\ena
All ${\mathbb C}$-linear relations among $(I;J)$'s are generated by
these relations. More precisely 
\vskip3mm

\begin{prop}\label{linear-rel}
$\ker\,(\,\tilde{ev}'\,)=\omega\wedge^{g-2} V$.
\end{prop}
\noindent
{\it Proof.} The proposition follows from Lemma \ref{injective} and 
Corollary \ref{resol-A}. \hfill \qed
\vskip2mm

The case $r=1$ of (\ref{lin-rel}) gives the obvious relations 
$\zeta_{ij}=\zeta_{ji}$.
The case $r=2$ occurs for $g\geq 4$ and gives 
\bea
&&
(i_1,j_1;j_2,j_3)-(i_1,j_2;j_1,j_3)+(i_1,j_3;j_1,j_2)=0.
\non
\ena
Corollary \ref{resol-A} tells that all ${\cal D}$-linear relations
among derivatives of $(I;J)$'s are obtained by differentiating
the relations
\bea
&&
\sum_{i=1}^{r+1}(-1)^{k-1}
\partial_i
(i_1,\ldots,i_r;j_1,\ldots,\hat{j_k},\ldots,j_{r+1})=0,
\label{d-rel}
\ena
and taking linear combinations of them.
The case $r=1$ of (\ref{d-rel}) gives the obvious relations
\bea
&&
\partial_k \zeta_{ij}=\partial_j\zeta_{ik},
\quad
j\neq k.
\label{d-rel-1}
\ena
The case $r=2$ gives the relations of the form
\bea
&&
\partial_m(ij;kl)-\partial_l(ij;km)+\partial_k(ij;lm)=0.
\label{d-rel-2}
\ena

In the following example we set 
\bea
&&
(I;J)_{i_1...i_p}=\partial_{i_1}\cdots\partial_{i_p}(I;J),
\qquad
(I;J)_{1^{r_1}...n^{r_n}}=\partial_1^{r_1}\cdots\partial_n^{r_n}(I;J),
\non
\\
&&
\zeta_{i_1...i_p}=\partial_{i_1}\cdots\partial_{i_p},
\log\,\theta(z),
\qquad
\zeta_{1^{r_1}...n^{r_n}}=
\partial_1^{r_1}\cdots\partial_n^{r_n}\log\,\theta(z).
\non
\ena

\vskip7mm
\noindent
{\bf \large Example} The case of $g=2$.
\vskip5mm

In this case $H^2\simeq W^2$ and 
\bea
&&
H^2=H^2_{-2}\oplus H^2_0\oplus H^2_1,
\non
\\
&&
H^2_{-2}=\mathbb{C}dz^2,
\quad
H^2_{0}=\mathbb{C}\zeta_{11}dz^2\oplus \mathbb{C}\zeta_{12}dz^2 \oplus 
\mathbb{C}\zeta_{22}dz^2,
\quad
H^2_1=\mathbb{C}(12;12)dz^2.
\non
\ena
By Theorem \ref{resol-grA} $\gr\,A$ is generated by $H^2(dz^2)^{-1}$
over ${\cal D}$. This means that
\bea
&&
\gr_0\,A={\mathbb C}1,
\quad
\gr_1\,A=\{0\},
\quad
\gr_2\,A=\mathbb{C}\zeta_{11}\oplus \mathbb{C}\zeta_{12} \oplus 
\mathbb{C}\zeta_{22},
\non
\ena
and that $\gr_n\, A$ for $n\geq 3$ is linearly spanned
by $\zeta_{i_1...i_{n}}$ and $(12;12)_{j_1...j_{n-3}}$.
More precisely the following elements form a $\mathbb{C}$-linear
basis of $gr_n\,A$:
\bea
&&
\zeta_{1^r2^s}\quad (r+s=n),
\qquad
(12;12)_{1^{r'}2^{s'}}
\quad (r'+s'=n-3).
\non
\ena
In fact $\gr_n\,A$ is the ${\mathbb C}$-span of these functions
as we said above. On the other hand the number of those functions
is $2n-1=\dim\,\gr_n\,A$. Thus they form a basis.
This basis for $n\leq 3$ is previously given in \cite{G}.
\vskip8mm

\noindent
{\bf \large Example} The case $g=3$.
\vskip5mm

We introduce the lexicographical order on the set ${\mathbb Z}^2$.
Then
\bea
&&
H^3\simeq W^3\oplus \mathbb{C} v dz^{3},\quad v\in \gr_2\,A,
\non
\\
&&
H^3_{-3}=\mathbb{C}dz^3,
\quad
H^3_{-1}=\oplus_{1\leq i,j\leq 3}\mathbb{C}\zeta_{ij}dz^3
\oplus \mathbb{C}vdz^3,
\non
\\
&&
H^3_0=\oplus_{1\leq i\leq j\leq 3,1\leq k\leq l\leq 3,(ij)<(kl)}
\, \mathbb{C}(ij;kl)dz^3,
\quad
H^3_1=\mathbb{C}(123;123)dz^3.
\non
\ena
Here we take a basis $\{v\}$ of $U^3$, $U^3=\mathbb{C} v dz^{3}$.
The relations (\ref{d-rel}) with $r=2$ in this case are
\bea
&&
\partial_1(ij;23)-\partial_2(ij;13)+\partial_3(ij;12)=0,
\quad
(ij)=(12),(13),(23).
\non
\ena
Using these relations functions of the form
$P(\partial)\partial_3(12;ij)$ can be eliminated 
from the ${\cal D}$-linear span of 
$H^3(dz^3)^{-1}$.
Then we get  $\mathbb{C}$-linear basis of $\gr_n\,A$ as
\bea
&&
\zeta_{1^{i_1}2^{i_2}3^{i_3}}\quad (i_1+i_2+i_3=n),
\qquad
v_{1^{i_1}2^{i_2}3^{i_3}}\quad (i_1+i_2+i_3=n-2),
\non
\\
&&
(12;12)_{1^{i_1}2^{i_2}},
\quad
(12;13)_{1^{i_1}2^{i_2}},
\quad
(12;23)_{1^{i_1}2^{i_2}},
\quad
(i_1+i_2=n-3),
\non
\\
&&
(13;13)_{1^{i_1}2^{i_2}3^{i_3}}
\quad
(13;23)_{1^{i_1}2^{i_2}3^{i_3}}
\quad
(23;23)_{1^{i_1}2^{i_2}3^{i_3}}
\quad
(i_1+i_2+i_3=n-3),
\non
\\
&&
(123;123)_{1^{i_1}2^{i_2}3^{i_3}}
\quad
(i_1+i_2+i_3=n-4).
\non
\ena
This can proved by checking that the number of the above functions
is $n^3-(n-1)^3=\dim\, \gr_n\,A$.

It is possible to give an explicit basis of $U^3$.
To this end we take Klein's sigma function 
$\sigma(u)=\sigma(u;\lambda)$, $u=(u_1,u_2,u_5)$ of the genus three non-hyperelliptic
curve $y^3-x^4-\sum_{0\leq \alpha<3,0\leq \beta<2}\lambda_{3\alpha+4\beta}
x^\alpha y^\beta=0$ (\cite{BEL1,BEL2}) as a theta function. We mainly follow the notations 
in \cite{BEL2} except indices of variables $u_i$. In this case we define
$\Theta=\left(\sigma(u)=0\right)$, $\partial_i=\partial/\partial u_i$
$\zeta_{i_1...i_k}=\partial_{i_1}\cdots\partial_{i_k}\log\,\sigma(u)$ for
$i_j\in\{1,2,5\}$, $(I;J)=\det(\zeta_{i_pj_q})_{1\leq p,q\leq r}$ 
for $I,J\in \{1,2,5\}^r$. The affine ring  $A$, in this case, is the space
of meromorphic functions on ${\mathbb C}^g$ which have poles only on
the zero set of $\sigma$ and are periodic with respect to the lattice
$2\omega_1{\mathbb Z}^3+2\omega_2{\mathbb Z}^3$ associated with $\sigma$.
All arguments in this paper do not
depend on a special choice of theta functions and linear coordinates.
Thus Theorem \ref{resol-grA} holds without any change.
Consequently a $\mathbb{C}$-linear basis of $\gr_n\,A$ is given by
\bea
&&
\zeta_{1^{i_1}2^{i_2}5^{i_5}}\quad (i_1+i_2+i_5=n),
\qquad
w_{1^{i_1}2^{i_2}5^{i_5}}\quad (i_1+i_2+i_5=n-2),
\non
\\
&&
(12;12)_{1^{i_1}2^{i_2}},
\quad
(12;15)_{1^{i_1}2^{i_2}},
\quad
(12;25)_{1^{i_1}2^{i_2}},
\quad
(i_1+i_2=n-3),
\non
\\
&&
(15;15)_{1^{i_1}2^{i_2}5^{i_5}}
\quad
(15;25)_{1^{i_1}2^{i_2}5^{i_5}}
\quad
(25;25)_{1^{i_1}2^{i_2}5^{i_5}}
\quad
(i_1+i_2+i_5=n-3),
\non
\\
&&
(125;125)_{1^{i_1}2^{i_2}5^{i_5}}
\quad
(i_1+i_2+i_5=n-4).
\non
\ena
Here one can take
\bea
&&
w=\frac{D_2^4\sigma\cdot \sigma}{\sigma(u;\lambda)^2}
=2(\zeta_{2222}+6\zeta_{22}^2),
\label{basis-u}
\ena
where $D_2^4$ is the Hirota derivative defined by
\bea
&&
D_2^4\sigma\cdot \sigma=
\frac{\partial^4}{\partial y_2^4}
\sigma(u+y;\lambda)\sigma(u-y;\lambda)\vert_{y=0}.
\non
\ena
This can be easily checked by using the expansion
\bea
&&
\sigma(u;\lambda)=S(u)+\sum_{n\geq 6} S_n(u),
\label{sigma-expansion}
\ena
where 
\bea
&&
S(u)=\frac{1}{20}(u_1^5-5u_1u_2^2+4u_5)
\non
\ena
and $S_n(u)$ is a polynomial in $u_1,u_2,u_5$ of degree $n$
with the degree being defined by $\deg\, u_i=i$.

\vskip2mm
\noindent
{\bf Remark} It is asserted that (\ref{sigma-expansion}) holds in 
\cite{BEL2,BEL3}.
However the proof of \cite{BEL3} contains some mistake. We could not find a proof of it
in a literature.  Thus, precisely speaking, the expression (\ref{basis-u})
is derived assuming the expansion (\ref{sigma-expansion}).
\vskip10mm

\noindent
{\bf Acknowledgements} This research is supported by Grant-in-Aid for Scientific
Research (B) 17340048. The results of this paper were presented at the work shop Integrable Systems, Geometry and Abelian Functions held at Tokyo Metropolitan
University, May 2005. We would like to thank Martin Guest, Yoshihiro Onishi and
Shigeki Matsutani for organizing the stimulating work shop.

\appendix
\section{Proof of Lemma \ref{sublem3}}
Condider the exact sequence
\bea
&&
0\lar {\cal O} \lar {\cal O}(1) \lar \gr_1\, {\cal O} \lar 0.
\label{a1}
\ena
By Lemma \ref{vanishing-lem} its cohomology sequence gives the isomorphism
\bea
&&
\beta: H^0\left(J,\gr_1\,{\cal O}\right)
\simeq
H^1\left(J,{\cal O}\right)
\simeq 
\bt^\ast.
\label{a2}
\ena

\begin{lemma}\label{a-lem1}
(i) The function $\zeta_i(z)$, $1\leq i\leq g$,
 can be naturally considered as an
element of $H^0\left(J,\gr_1\,{\cal O}\right)$ and 
form its linear basis, where $\zeta_i(z)$ is defined by (\ref{zeta-i}).

\noindent
(ii) $\beta\left(\zeta_j(z)\right)=
\pi \sum_{k=1}^g \left((Im\,\tau)^{-1}\right)_{jk}d{\bz}_k.$
\end{lemma}
\vskip2mm
\noindent
{\it Proof.}
Let $p:V\lar V/\Gamma=J$ be the projection, where 
$\Gamma={\mathbb Z}^g+\tau {\mathbb Z}^g$.
We shall use the group cohomology description of the cohomology 
groups of sheaves on
$J$ by pulling them back to $V$. We refer the appendix to \S 2 
in the book \cite{M}. We mainly follow the notations in that book.

For a sheaf ${\cal F}$ on J such that $H^i(V,p^\ast {\cal F})=0$, $i\geq 1$,
there is an isomorphism
\bea
&&
H^p\left(\Gamma, H^0(V,p^\ast {\cal F})\right)
\simeq
H^p(J,{\cal F}),
\quad
p\geq 0,
\label{a3}
\ena
where the left hand side is the group cohomology with the
value in the $\Gamma$-module $H^0(V,p^\ast {\cal F})$. 
Since $V$ is a Stein manifold and $p$ is a local isomorphism,
 (\ref{a3}) is valid for a coherent ${\cal O}_J$-module ${\cal F}$.
Applying (\ref{a3}) we have
\bea
&&
H^0(J,\gr_1\,{\cal O})
\simeq
H^0\left(\Gamma,H^0(V,p^\ast \gr_1\,{\cal O})\right)
=H^0(V,p^\ast \gr_1\,{\cal O})^\Gamma,
\label{a4}
\ena
where $(\quad)^\Gamma$ denotes the set of $\Gamma$-invariants.

We pull back the sequence (\ref{a1}) to $V$ and take the cohomology sequence
of it. Then we get
\bea
&&
0 \lar H^0(V,{\cal O}_V)
\lar
H^0\left(V, {\cal O}(p^\ast\Theta)\right)
\lar
H^0\left(V, p^\ast \gr_1\, {\cal O}\right)
\lar
0,
\label{a41}
\ena
since $H^i(V,{\cal O}_V)=0$, $i\geq 1$.
Thus
\bea
&&
H^0\left(V, p^\ast \gr_1\, {\cal O}\right)
\simeq
\frac{H^0\left(V, {\cal O}(p^\ast\Theta)\right)}{H^0(V,{\cal O}_V)}.
\non
\ena
Then
\bea
&&
H^0(J,\gr_1\,{\cal O})
\simeq
\left(
\frac{H^0\left(V, {\cal O}(p^\ast\Theta)\right)}{H^0(V,{\cal O}_V)}
\right)^\Gamma.
\label{a5}
\ena

Since $\zeta_j(z)\in H^0\left(V,{\cal O}_V(p^\ast \Theta)\right)$ and
\bea
&&
\zeta_j(z+m+\tau n)=\zeta_j(z)-2\pi i n_j,
\quad
m, n\in \mathbb{Z}^g,
\non
\ena
$\zeta_j(z)$ can be considered as an element of the right hand side of
(\ref{a5}) and therefore of $H^0(J,\gr_1\,{\cal O})$. We define the element $a_j(m+\tau n,z)$ in the one cocycle
$C^1\left(\Gamma,H^0(V,{\cal O}_V)\right)$ by
\bea
&&
a_j(m+\tau n,z)=-2\pi i n_j.
\label{aj}
\ena
Let $\delta_1$ be the connecting homomorphism
\bea
\delta_1: 
H^0\left(\Gamma,H^0(V,p^\ast \gr_1\,{\cal O})\right)
\lar
H^1\left(\Gamma,H^0(V,{\cal O}_V)\right)
\non
\ena
of the group cohomology sequence of (\ref{a41}).
One can easily check that
\bea
&&
\delta_1\left(\zeta_j(z)\right)=a_j(m+\tau n,z).
\label{a6}
\ena
Let us calculate the description of $a_j(m+\tau n,z)$ in terms of
Dolbeaut cohomology description.

Let ${\cal C}^{p,q}_X$ be the sheaf of germs of $C^\infty$ $(p,q)$-forms on
$X$.
Notice that
\bea
&&
H^1\left(\Gamma,H^0(V,{\cal O}_V)\right)
\simeq
H^1(J,{\cal O}_J)
\simeq
\frac{\left(\ker\left(\bar{\partial}: 
H^0(V,{\cal C}^{0,1}) \lar H^0(V,{\cal C}^{0,2}_V)\right)\right)^\Gamma}
{\bar{\partial}\,\left( H^0(V,{\cal C}^{0,0}_V)^\Gamma\right)}.
\label{a7}
\ena
This isomorphism is given in the following way.

We have a natural map
\bea
&&
H^1\left(\Gamma,H^0(V,{\cal O}_V)\right)
\lar
H^1\left(\Gamma,H^0(V,{\cal C}_V^{0,0})\right)
\simeq
H^1(J,{\cal C}^{0.0}_J)
=0.
\non
\ena
This means that any one cocycle 
$f(\gamma,z)\in C^1\left(\Gamma,H^0(V,{\cal O}_V)\right)$ is written
as
\bea
&&
f(\gamma,z)=h(z+\gamma)-h(z),
\non
\ena
for some $C^\infty$ function $h(z)$ on $V$.
One can easily check that the isomorphism (\ref{a7}) is given
by
\bea
&&
f(\gamma,z) \mapsto \bar{\partial} h(z).
\non
\ena
For the one cocycle $a_j(\gamma,z)$ we have
\bea
&&
a_j(\gamma,z)=h_j(z+\gamma)-h_j(z),
\quad
\gamma\in \Gamma,
\non
\\
&&
h_j(z)=\pi \sum_{k=1}^g \left((Im\,\tau)^{-1}\right)_{jk}(\bz_k-z_k).
\non
\ena
Thus the isomorphism (\ref{a7}) is given by
\bea
&&
a_j(\gamma,z)
\mapsto
\pi \sum_{k=1}^g \left((Im\,\tau)^{-1}\right)_{jk}d\bz_k.
\label{a71}
\ena
This proves (ii) of Lemma \ref{a-lem1}. Since $\beta$ is an isomorphism
and $\beta\left(\zeta_j(z)\right)$, $1\leq j\leq g$ are linearly independent
due to (\ref{a71}), (i) of Lemma \ref{a-lem1} is also proved.
\hfill\qed
\vskip2mm

By Lemma \ref{a-lem1} and (\ref{free2}) we have
\bea
&&
\Xi=d\log\, \theta =\sum_{j=1}^g\zeta_j(z) dz_j
\in
H^0(J,\gr_1\,{\cal O})\otimes H^0(J,\Omega^1)\simeq H^0(J,\gr_1\,\Omega^1).
\non
\ena
Consider the connecting homomorphism
\bea
\tilde{\beta}:H^0(J,\gr_1\,\Omega^1) \simeq 
H^1(J,\Omega^1)\simeq \bt^\ast\otimes T^\ast.
\label{a8}
\ena
of the cohomology sequence of
\bea
&&
0
\lar
\Omega^1
\lar
\Omega^1(1)
\lar
\gr_1\, \Omega^1
\lar
0.
\non
\ena
By Lemma \ref{a-lem1} we have
\bea
&&
\tilde{\beta}(\Xi)
=\sum_{j=1}^g\beta\left(\zeta_j(z)\right)\wedge dz_j
=
\pi
\sum_{i,j=1}^g \left((Im\, \tau)^{-1}\right)_{ij}d\bz_i\wedge dz_j,
=\hat{\omega}
\label{a9}
\ena
where we use the fact that $\tau$ is symmetric. 

Consider the composition of maps
\bea
&&
\Omega^{g-n}
\lar
\gr_0\, \Omega^{g-n}
\stackrel{\Xi\wedge}{\lar}
\Xi\wedge \gr_0\,\Omega^{g-n}
\hookrightarrow
\gr_1\,\Omega^{g-n+1},
\non
\ena
where the first map is the natural projection.
It induces maps of cohomologies
\bea
&&
H^i(J,\Omega^{g-n})
\simeq
H^i(J,\gr_0\,\Omega^{g-n})
\simeq
H^i(J,\Xi \wedge \gr_0\,\Omega^{g-n})
\lar
H^i(J,\gr_1\,\Omega^{g-n+1}),
\quad
i\leq n-2,
\label{a10}
\ena
where we use Lemma \ref{fundamental-cohomology} and Lemma \ref{sublem2}.
Let us study the map (\ref{a10}).
In sum the map (\ref{a10}) is induced from the map
\bea
&&
\Omega^{g-n}
\stackrel{\Xi\wedge}{\lar}
\gr_1\,\Omega^{g-n+1}.
\label{a11}
\ena
Thus we consider more generally the bilinear map
\bea
&&
\Omega^{g-n}\times \gr_1\,\Omega^1 \lar \gr_1\,\Omega^{g-n+1}.
\label{a12}
\ena
We pull it back to $V$ and get the map of $\Gamma$-modules
\bea
&&
H^0(V,\Omega_V^{g-n})\times H^0(V,p^\ast \gr_1\,\Omega^1)
\lar
H^0(V,p^\ast\gr_1\,\Omega^{g-n+1}).
\label{a13}
\ena
It induces a cup product on cohomologies (\cite{M})
\bea
&&
H^i\left(\Gamma,H^0(V,\Omega_V^{g-n})\right)
\times
H^j\left(\Gamma,H^0(V,p^\ast \gr_1\,\Omega^1)\right)
\lar
H^{i+j}\left(\Gamma,H^0(V,p^\ast\gr_1\,\Omega^{g-n+1})\right).
\label{a14}
\ena
Taking $j=0$ we have a map
\bea
&&
H^i\left(\Gamma,H^0(V,\Omega_V^{g-n})\right)
\times
H^0\left(\Gamma,H^0(V,p^\ast \gr_1\,\Omega^1)\right)
\lar
H^{i}\left(\Gamma,H^0(V,p^\ast\gr_1\,\Omega^{g-n+1})\right).
\label{a15}
\ena
Then $\Xi\in H^0\left(\Gamma,H^0(V,p^\ast \gr_1\,\Omega^1)\right)$
defines a map
\bea
&&
H^i\left(\Gamma,H^0(V,\Omega_V^{g-n})\right)
\lar
H^{i}\left(\Gamma,H^0(V,p^\ast\gr_1\,\Omega^{g-n+1})\right),
\label{a16}
\ena
which, by the definition, coincides with the map (\ref{a10}).
Let us describe this map in terms of Dolbeaut cohomology.
To this end consider the diagram
\bea
&&
\begin{array}{ccccc}
H^i(J,\Omega^{g-n}) & \times & H^0(J,\gr_1\,\Omega^1) & \lar &
H^i(J,\gr_1\,\Omega^{g-n+1})
\\
\quad & \mapdown{1\times \tilde{\beta}} & \quad & \quad & \mapdown{\delta_2}
\\
H^i(J,\Omega^{g-n}) & \times & H^1(J,\Omega^1) & \lar &
H^{i+1}(J,\Omega^{g-n+1}).
\end{array}
\label{a17}
\ena
Here the down arrows are isomorphisms, $\delta_2$ is the connecting
homomorphism of the long cohomology exact sequence of
\bea
&&
0\lar \Omega^{g-n+1} \lar \Omega^{g-n+1}(1) \lar \gr_1\,\Omega^{g-n+1}
\lar 0,
\non
\ena
and the horizontal map in the second row is the cup product 
of cohomologies.
By a direct calculation in terms of group cohomologies one can show
that (\ref{a17}) is a commutative diagram.
The cup product of group cohomologies is compatible with
that of sheaf cohomologies. In Dolbeaut description the cup
product is given by the exterior product. Therefore, by (\ref{a9}), the map
$\alpha$ is given by wedging $\hat{\omega}$.
\hfill\qed

\section{Proof of (\ref{omega=omega})}

\begin{lemma}\label{b1}
As an element of $H^1(J,\mathbb{C})$ we have
\bea
&&
d\zeta_j=\pi \sum_{k=1}^g\left((Im\,\tau)^{-1}\right)_{jk}(d\bz_k-dz_k).
\label{b2}
\ena
\end{lemma}
\vskip2mm
\noindent
{\it Proof.}
We denote the one cycle on $J$ specified by the element $m+\tau n$ in
$\Gamma$ by $\gamma(m,n)$.
Then
\bea
&&
\int_{\gamma(m,n)}d\zeta_j=-2\pi i n_j,
\quad
\int_{\gamma(m,n)}(d\bz_k-dz_k)=-2i\sum_{j=1}^gIm\,\tau_{kj}n_j.
\non
\ena
The lemma follows from this.
\hfill\qed
\vskip2mm

By Lemma \ref{b1} we have
\bea
\omega=\sum_{j=1}^g d\zeta_j\wedge dz_j
&=&
\pi \sum_{j,k=1}^g\left((Im\,\tau)^{-1}\right)_{jk}(d\bz_k-dz_k)\wedge dz_j
\non
\\
&=&
\pi \sum_{j,k=1}^g\left((Im\,\tau)^{-1}\right)_{jk}d\bz_k \wedge dz_j
\non
\\
&=&
\hat{\omega},
\non
\ena
where we use the symmetry of $Im\, \tau$.
\hfill\qed

\section{Proof of Lemma \ref{injective}}
For $n=0,1,g+1$ the lemma is obvious.
We assume $2\leq n\leq g$.
By Lemma \ref{fundamental-cohomology}, Lemma \ref{sheaf-lemma} (iii)
and Lemma \ref{a-lem1} (i)
\bea
&&
\bt^\ast\simeq H^1(J,{\cal O})\stackrel{\delta_1^{-1}}{\simeq}
H^0(J,\gr_1\,{\cal O})
\stackrel{d}{\simeq}
H^0(J,d\gr_1\,{\cal O})= \oplus_{i=1}^g {\mathbb C}d\zeta_i\simeq V_{+},
\label{hw-4}
\ena
where $\delta_1$ is the connecting homomorphism of the cohomology sequence
of (\ref{a1}) and the last isomorphism 
is given by $d\zeta_j\mapsto \beta_j$.
On the other hand $T^\ast=\oplus_{i=1}^g\mathbb{C}dz_i\simeq V_{-}$ by the map
$dz_j\mapsto \alpha_j$. Thus we have, 
by (\ref{hi-1}) and (\ref{hi-2}),
\bea
H^i\left(J,\Xi\wedge\gr_0\,\Omega^{g-n}\right)
&\simeq&
 \wedge^i V_{+}\wedge^{g-n} V_{-}
\quad
i\leq n-2,
\label{hw-5}
\\
H^i\left(J,\gr_1\,\Omega^{g-n+1}\right)
&\simeq&
 \wedge^{i+1} V_{+}\wedge^{g-n+1} V_{-}
\quad
i\geq 0.
\label{hw-6}
\ena
By Lemma \ref{sublem3}, Lemma \ref{a-lem1} and (\ref{hw-4}) the map
\bea
&&
\alpha:
H^i\left(J,\Xi\wedge\gr_0\,\Omega^{g-n}\right)
\lar
H^i\left(J,\gr_1\,\Omega^{g-n+1}\right)
\non
\ena
is given, in the description of (\ref{hw-5}) and (\ref{hw-6}), by
wedging $\omega=\sum_{j=1}^g \beta_j \wedge \alpha_j$.

Consider the exact sequence (\ref{h-exact4}).
For $n=2$ we have the exact sequence
\bea
&&
0\lar \wedge^{g-2} V_{-} \stackrel{\wedge\omega}{\lar}
V_{+}\wedge^{g-1} V_{-} \lar
H^0(J,d\gr_1\,\Omega^{g-1}) \lar \cdots.
\label{hw-7}
\ena
Since
\bea
&&
d\gr_1\,\Omega^{g-1}\simeq \gr_2\,\Omega^g,
\quad
H^0(J,\gr_2\,\Omega^g)\simeq H^g_{2-g},
\label{hw-8}
\ena
we get the injective map
\bea
&&
W^g_{2-g}=\frac{V_{+}\wedge^{g-1} V_{-}}{\omega\wedge^{g-2} V_{-}}
\lar
H^g_{2-g}.
\label{hw-9}
\ena
This map obviously coincides with $\varphi_2$. Thus the $n=2$ case of
the lemma is proved.

Let us assume $n\geq 3$. By the isomorphisms (\ref{isom2}) with
$k=n-1$ and (\ref{isom-1}) we have
\bea
&&
H^{n-2}(J,d\gr_1\,\Omega^{g-n+1})
\simeq
H^1(J,d\gr_{n-2}\,\Omega^{g-2})
\simeq
H^g_{n-g},
\quad
n\geq 3.
\label{hw-10}
\ena
Thus we have the injective map
\bea
&&
\psi_n:W^g_{n-g} \lar H^g_{n-g}.
\label{hw-11}
\ena
Let us prove that $\psi_n$ coincides with $\varphi_n$.
To this end we explicitly describe the isomorphism (\ref{hw-10}).
We use the notations in Appendix A.

\begin{lemma}\label{c-lem1}
Let $I$, $J$ be the index sets such that $|I|=n\geq 1$, $|J|=l$,
$n+l\leq g$. Write $I=(i_1,\ldots,i_n)$. Then
\bea
&&
d\zeta_I\wedge dz_J=d(\eta\wedge dz_J),
\label{hw-12}
\\
&&
\eta=\frac{1}{n}\sum_{k=1}^n
(-1)^{k-1}
\zeta_{i_k}
d\zeta_{i_1}\wedge\cdots\wedge \widehat{d\zeta_{i_k}}\wedge 
\cdots\wedge d\zeta_{i_n},
\non
\ena
and $\eta\in H^0\left(V,\Omega_V^{n-1}(np^\ast \Theta)\right)$.
In particular 
\bea
&&
d\zeta_I\wedge dz_J\in H^0\left(V,d\Omega_V^{n+l-1}(np^\ast \Theta)\right).
\non
\ena
\end{lemma}
\vskip2mm
\noindent
{\it Proof.}
Notice that the determinant
$$
\left|
\begin{array}{cccc}
\zeta_{i_1}&\zeta_{i_1r_1}&\cdots &\zeta_{i_1r_{n-1}}\\
\vdots&\vdots&\quad&\vdots\\
\zeta_{i_n}&\zeta_{i_nr_1}&\cdots &\zeta_{i_nr_{n-1}}\\
\end{array}
\right|
\non
$$
has poles on $p^\ast\Theta$ of order at most $n$ for any 
$(r_1,\ldots,r_{n-1})$. The lemma follows from this.
\hfill\qed
\vskip2mm

Let $a_j(\gamma)$ be the one cocycle defined in (\ref{aj}) and
\bea
&&
f_{IJ}(\gamma_1,\ldots,\gamma_k|z)
=a_{r_1}(\gamma_1)\cdots a_{r_k}(\gamma_k)d\zeta_I\wedge dz_J,
\label{hw-14}
\ena
where $r_1,\ldots,r_k$ are distinct. We omit writing $r_j's$ in $f_{IJ}$
for the sake of simplicity.
By Lemma \ref{c-lem1}, for $I$, $J$ satisfying $|I|=n-k-1$, $|J|=g-n+1$,
$f_{IJ}(\gamma_1,\ldots,\gamma_k|z)$ can be considered as an element of
\bea
&&
H^k\left(\Gamma,H^0(V,p^\ast d\gr_{n-k-1}\,\Omega^{g-k-1})\right)
\simeq
H^k(J,d\gr_{n-k-1}\,\Omega^{g-k-1}).
\label{hw-15}
\ena
Let $\iota_k$ be the connecting homomorphism (\ref{isom1})
\bea
&&
\iota_k:
H^k\left(\Gamma,H^0(V,p^\ast d\gr_{n-k-1}\,\Omega^{g-k-1})\right)
\simeq
H^{k+1}\left(\Gamma,H^0(V,p^\ast d\gr_{n-k-2}\,\Omega^{g-k-2})\right).
\label{hw-16}
\ena
Then, using Lemma \ref{c-lem1}, we have
\bea
&&
\iota_k(f_{IJ})(\gamma_0,\ldots,\gamma_k)
=
\non
\\
&&
a_{r_1}(\gamma_1)\cdots a_{r_k}(\gamma_k)
\frac{1}{n-k-1}
\sum_{l=1}^{n-k-1}
(-1)^{l-1}
a_{i_l}(\gamma_0)
d\zeta_{i_1}\wedge\cdots\wedge \widehat{d\zeta_{i_l}}\wedge 
\cdots\wedge d\zeta_{i_{n-k-1}}\wedge dz_J.
\label{hw-17}
\ena

With the help of (\ref{hw-17}) we find that the composition of maps
\bea
&&
H^0(J,\gr_n\Omega^g)\simeq
H^0(J,d\gr_{n-1}\Omega^{g-1})\lar
H^1(J,d\gr_{n-2}\Omega^{g-2})\simeq
H^{n-2}(J,d\gr_{1}\Omega^{g-n+1})
\non
\ena
is given by
\bea
&&
d\zeta_I\wedge dz_J \mapsto
\frac{1}{n-1}\sum_{k=1}^{n-1}
(-1)^{n-1-k}a_{I\backslash\{i_k\}}(\gamma_{n-2},...,\gamma_1)
d\zeta_{i_k}\wedge dz_J,
\label{half-1}
\ena
where $I=(i_1,...,i_{n-1})$, $|J|=g-n+1$, and
\bea
&&
a_{I\backslash\{i_k\}}(\gamma_{n-2},...,\gamma_1)
=a_{i_1}(\gamma_1)\cdots a_{i_{k-1}}(\gamma_{k-1})a_{i_{k+1}}(\gamma_k)
\cdots a_{i_{n-1}}(\gamma_{n-2}).
\non
\ena
On the other hand the map
\bea
&&
\wedge^{n-1} H^1(J,{\cal O})\otimes H^0(J,\Omega^{g-n+1})
\simeq
H^{n-2}(J,\gr\Omega^{g-n+1})
\non
\ena
is given by
\bea
&&
\zeta_Idz_J\mapsto 
\frac{1}{n-1}\sum_{k=1}^{n-1}
(-1)^{n-1-k}a_{I\backslash\{i_k\}}(\gamma_{n-2},...,\gamma_1)
\zeta_{i_k}\wedge dz_J
\non
\ena
and consequently the composition of maps
\bea
&&
\wedge^{n-1} V_{+}\wedge^{g-n+1} V_{-}
\simeq
\wedge^{n-1} H^1(J,{\cal O})\otimes H^0(J,\Omega^{g-n+1})
\simeq
H^{n-2}(J,\gr\Omega^{g-n+1})
\lar
H^{n-2}(J,d\gr_{1}\Omega^{g-n+1})
\non
\ena
is given by
\bea
&&
\beta_I\wedge \alpha_J
\mapsto
\frac{1}{n-1}\sum_{k=1}^{n-1}
(-1)^{n-1-k}a_{I\backslash\{i_k\}}(\gamma_{n-2},...,\gamma_1)
d\zeta_{i_k}\wedge dz_J.
\label{half-2}
\ena
Comparing (\ref{half-1}) with (\ref{half-2}) we have $\psi_n=\varphi_n$.
\hfill\qed

\end{document}